\documentclass[a4paper,12pt]{article}
    \usepackage[top=2.5cm,bottom=2.5cm,left=2.5cm,right=2.5cm]{geometry}
    \usepackage{cite, amsmath, amssymb}
   \usepackage{tikz}
\begin{document}
$\vspace{3cm}$
\begin{center}
{\LARGE\bf $L$-Borderenergetic Graphs and Normalized Laplacian Energy}
\end{center}
\begin{center}
{\large \bf Fernando Tura}
\end{center}
\begin{center}
\it Departamento de Matem\'atica, UFSM, Santa Maria, RS, 97105-900, Brazil
\end{center}
\begin{center}
\tt ftura@smail.ufsm.br
\end{center}



\newtheorem{Thr}{Theorem}
\newtheorem{Pro}{Proposition}
\newtheorem{Que}{Question}
\newtheorem{Con}{Conjecture}
\newtheorem{Cor}{Corollary}
\newtheorem{Lem}{Lemma}
\newtheorem{Fac}{Fact}
\newtheorem{Ex}{Example}
\newtheorem{Def}{Definition}
\newtheorem{Prop}{Proposition}
\def\floor#1{\left\lfloor{#1}\right\rfloor}

\newenvironment{my_enumerate}{
\begin{enumerate}
  \setlength{\baselineskip}{14pt}
  \setlength{\parskip}{0pt}
  \setlength{\parsep}{0pt}}{\end{enumerate}
}

\newenvironment{my_description}{
\begin{description}
  \setlength{\baselineskip}{14pt}
  \setlength{\parskip}{0pt}
  \setlength{\parsep}{0pt}}{\end{description}
}


\begin{abstract}
The Laplacian and normalized Laplacian energy of $G$ are given by  expressions $E_{L}(G) = \sum_{i=1}^n | \mu_i - \overline{d}|,   E_{\mathcal{L}}(G)= \sum_{i=1}^n | \lambda_i - 1 |,$ respectively, where   $\mu_i$ and $\lambda_i$  are the eigenvalues of Laplacian matrix $L$ and normalized Laplacian matrix $\mathcal{L}$ of $G.$ An interesting problem in {\em spectral graph theory} is to find graphs $\{L, \mathcal{L} \}-$noncospectral with the same   $E_{ \{L, \mathcal{L}\}}(G).$
In this paper, we  present graphs of order $n,$ which are $L$-borderenergetic (in short, $E_{L}(G) = 2n-2)$
and graphs  $\mathcal{L}$-noncospectral  with the same normalized Laplacian energy.
\end{abstract}

\baselineskip=0.30in

\begin{center}
     (Received November 4, 2016)
\end{center}

\section{Introduction}
\label{intro}
Throughout this paper, all graphs are assumed to be finite, undirected and without loops or multiple edges.
If $G$ is a graph of order $n$ and $M$ is a real symmetric matrix associated with $G,$ then  the $M$- {\em energy} of $G$ is
\begin{equation}
\label{energy}
 E_{M}(G)= \sum_{i=1}^n |  \lambda_i (M) - \frac{ tr(M)}{n}|.
 \end{equation}
The  energy $E(G)$ of a graph $G$ simply refers to using  the adjacency matrix in (\ref{energy}).
There are many results on energy and its applications in several areas, including in chemistral
 see  \cite{Gutman2012} for more details and the references \cite{Li, Li2, Li3, Li4, Li5, Li6}.

Recently, a new concept as {\em borderenergetic}  graphs  \cite{Gutman2015} was proposed, namely graphs of order $n$ satisfying   $E(G) = 2n-2.$ In this way, several authors have been presented families of borderenergetic graphs \cite{JTT2015, Li, Li2, Hou, Li3}.

An analogous concept as borderenergetic graphs, called $L$-borderenergetic graphs was proposed   
in \cite{tura}. That is, a graph $G$ of order $n$ is $L$-borderenergetic if $E_{L}(G) =2n-2,$ where  
$E_{L}(G) = \sum_{i=1}^n | \mu_i - \overline{d}|,$ and $\overline{d}$ is the avarage degree of $G.$ Some  classes of $L$-borderenergetic of order $n=4r+4$ $(r\geq 1)$ are obtained in \cite{tura}. In \cite{Li7}, a kind of threshold graphs were found to be $L$-borderenergetic, and all the connected  non-complete and pairwise non-isomorphic $L$-borderenergetic graphs of small order $n$ depicted for $n$ with  $4\leq n \leq 9.$  

Since that finding noncospectral graphs with the same energy is an interesting problem in spectral graph theory, in this paper we continue this investigation  presenting some new graphs which are $L$- borderenergetic and finish it showing two classes of graphs that are $\mathcal{L}$-noncospectral  with the same normalized Laplacian energy.

The paper is organized as follows.
In Section 2 we describe some  known results about the Laplacian and normalized Laplacain spectrum of graphs.
In Section 3 we present two classes of  $L$-borderenergetic. We finalize this paper, showing graphs with the same normalized Laplacian energy.

\section{Premilinares}
\label{Diag}

Let $G_1= (V_1, E_1)$ and $G_2= (V_2, E_2)$  be undirected graphs without  loops or multiple edges.
The {\em union} $G_1 \cup G_2$ of graphs $G_1$ and $G_2$  is the graph $G=(V,E)$  for which $V= V_1 \cup V_2$ and $E = E_1 \cup E_2.$
We denote the graph  $ \underbrace{G  \cup G \cup \ldots \cup G}_{m} $ by  $m G.$
The {\em join} $G_1 \nabla G_2$  of graphs $G_1$ and $G_2$  is the graph obtained from $G_1 \cup G_2$ by joining every vertex of $G_1$
with every vertex of $G_2.$

The Laplacian spectrum of $G_1  \cup \ldots \cup G_k$ is the union of Laplacian spectra of $G_1, \ldots, G_k,$ while
the Laplacian spectra of the complement  of $n$- vertex graph $G$ consists  of values $n - \mu_i,$ for each Laplacian eigenvalue
$\mu_i$ of $G,$ except for a single instance of eigenvalue $0$ of $G.$

\begin{Thr}
\label{theorem2} (\cite{Merris2})
Let $G_1$ and $G_2$ be graphs on  $n_1$ and $n_2$ vertices, respectively. Let $L_1$ and $L_2$ be the Laplacian matrices for $G_1$ and $G_2,$
respectively, and let $L$ be the Laplacian matrix for $G_1 \nabla G_2.$  If $0=\alpha_1 \leq \alpha_2 \leq \ldots \leq \alpha_{n_1} $ and $0=\beta_1 \leq \beta_2 \leq \ldots \leq \beta_{n_2}$
are the eigenvalues of $L_1$ and $L_2,$ respectively. Then the eigenvalues of $L$ are
$$ 0, \hspace{0,2cm} n_2 + \alpha_2, \hspace{0,2cm} n_2 + \alpha_3 , \ldots, \hspace{0,1cm} n_2 + \alpha_{n_1}$$
$$ n_1 + \beta_2, \hspace{0,2cm} n_1 + \beta_3, \ldots, \hspace{0,1cm} n_1 + \beta_{n_2}, \hspace{0,2cm} n_1 + n_2.$$
\end{Thr}

The following result is due to Butler( \cite{Butler}, Theorem 12).

\begin{Thr}
\label{main4}
Let  $G_1 = (V_1, E_1)$ be an $r$- regular graph on $n$ vertices and $G_2 = (V_2, E_2)$ be an $s$- regular graph on $m$ vertices. Suppose 
$$ 0 = \lambda_1 \leq \ldots \leq \lambda_n \leq 2$$
are the $\mathcal{L}$- eigenvalues of $G_1$ and 
$$ 0 = \mu_1 \leq \ldots \leq \mu_m \leq 2$$
are the $\mathcal{L}$- eigenvalues of $G_2.$ Then the eigenvalues  of $G_1\nabla G_2$  are
$$ 0, \frac{m+r\lambda_2}{m+r}, \ldots, \frac{m+r\lambda_n}{m+r}, \frac{n+s\mu_2}{n+s}, \ldots, \frac{n+s\mu_m}{n+s}, \frac{m}{m+r}  + \frac{n}{n+s}.$$
\end{Thr}

\section{Constructing new $L$-Borderenergetic graphs}
\label{SecInterval}
Recall that the {\em Laplacian energy}  $E_{L}(G)$ of $G$ is  defined to be
$ \sum_{i=1}^n | \mu_i - \overline{d} |,$
 where $ 0= \mu_1 \leq \mu_2 \leq \ldots \leq \mu_n$  are the Laplacian eigenvalues  of $G$ and $\overline{d}$ is the average degree of $G.$
 It is known that the complete graph $K_n$  has Laplacian energy  $2(n-1).$
We exhibit two infinite classes of graphs which are $L$-noncospectral and $L$-borderenergetic.  

\subsection{ The class $ K_{n-1}\odot K_{n}$ }
\begin{figure}[h!]
\begin{center}
\begin{tikzpicture}
   [scale=0.6,auto=left,every node/.style={circle,inner sep=2pt}]
  \foreach \i/\w in {11/,12/,13/,14/,15/,16/,17/,18/}{
    \node[draw,circle,fill=black]   (\i) at ({72*(\i -1)+90}:2) {}; }
  \foreach \from in {12,13,14,15,16,17,18}{
    \foreach \to in {11,12,...,\from}
      \draw (\from) -- (\to);}

    \node[draw,circle,fill=black,label=below:$$] (u) at (-6,0.5) {};
    \node[draw,circle,fill=black,label=below:$$] (q) at (-7,-1.5) {};
  \node[draw,circle,fill=black,label=below:$$] (r) at (-5,-1.5) {};

\path

    (u) edge node[right]{}(12)
       (q) edge node[right]{}(12)
   (r) edge node[right]{}(12)

  (u) edge node[below]{}(q)
 (u) edge node[below]{}(r)
 (q) edge node[right]{}(r);

\end{tikzpicture}
          \caption{Graph $K_4 \odot K_5$}
       \label{fig1}
\end{center}
\end{figure}
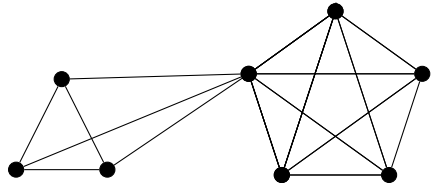

For  each integer $n\geq 3,$  we define the graph $G$  in  $K_{n-1} \odot K_{n}$ to be the following join
$$ G = (K_{n-1} \cup K_{n-2}) \nabla K_1$$
 of order $2n -2.$ Let $\mu^m$ denote  the laplacian eigenvalue $\mu $ with multiplicity equals to $m.$
The Figure \ref{fig1} shows the graph $K_4 \odot K_5.$

\begin{Lem}
\label{lema1}
Let $G =K_{n-1} \odot K_{n}$ be a graph of order $2n-2.$
Then  the Laplacian spectrum of $ G$ is given by
$$ 0 ;  \hspace{0,2cm} 1 ;\hspace{0,2cm} (n-1)^{n-3}; \hspace{0,2cm} n^{n-2} ;\hspace{0,2cm}  2n-2.$$

\end{Lem}
{\bf Proof:}
Let $G=K_{n-1} \odot K_{n}$ be a graph of order $2n-2.$ Since that $K_{n-1}$ and $K_{n-2}$ have Laplacian spectrum equal to $\{ (n-1)^{n-2},0\}$ and $\{ (n-2)^{n-3},0\},$ respectively. Taking $G_1= K_{n-1} \cup K_{n-2}$ and $G_2=K_1,$
according by Theorem \ref{theorem2}, follows that the Laplacian spectrum of $G$ is equal to
$$ 0 ;  \hspace{0,2cm} 1 ;\hspace{0,2cm} (n-1)^{n-3}; \hspace{0,2cm} n^{n-2} ;\hspace{0,2cm}  2n-2.$$

\begin{Thr}
\label{}
For each $n\geq 3,$  $G = K_{n-1} \odot K_{n-2}$ is  $L$-borderenergetic and  $L$-noncospectral graph with $K_{2n-2}.$
\end{Thr}
{\bf Proof:}
Clearly $G$ and $K_{2n-2}$ are  $L$-noncospectral.
Let $\overline{d}$ be the average degree of $G.$ Since that  $\overline{d}$ is equal to average of Laplacian eigenvalues of $G$ then
 $ \overline{d} = \frac{ 2n-2 + n(n-2) +(n-1)(n-3) +1}{2n-2} = n-1.$ Using Lemma \ref{lema1},
 $ E_{L}(G) = |2n -2 -(n-1)| + (n-2) | n -(n-1)| + (n-3) |(n-1)-(n-1)| +|1 -(n-1)| +|0-(n-1)| = 4n-6 = E_L(K_{2n-2}). $

\subsection{ The class $K_n \cdot K_n$ }
For  each integer $n\geq 3,$  we define the graph $G$ in  $K_n \cdot K_n$  of order $2n,$ obtained from two copies of the complete graph by adding $n$ edges between one vertex  of a copy of $K_n$ and $n$ vertices  of the other copy.

\begin{figure}[h!]
\begin{center}
\begin{tikzpicture}
    [scale=0.6,auto=left,every node/.style={circle,inner sep=2pt}]
     \node[draw,circle,fill=black,label=below:$$] (1) at (-2,0.5) {};  
    \node[draw,circle,fill=black,label=below:$$] (2) at (-3,-1.5) {};  
     \node[draw,circle,fill=black,label=below:$$] (3) at (-1,-1.5) {};

    \node[draw,circle,fill=black,label=below:$$] (u) at (-6,0.5) {};
    \node[draw,circle,fill=black,label=below:$$] (q) at (-7,-1.5) {};
  \node[draw,circle,fill=black,label=below:$$] (r) at (-5,-1.5) {};

\path

    (u) edge node[right]{}(1)
       (q) edge node[right]{}(1)
   (r) edge node[right]{}(1)

     (1) edge node[below]{}(2)
     (1) edge node[below]{}(3)
      (2) edge node[right]{}(3)

  (u) edge node[below]{}(q)
 (u) edge node[below]{}(r)
 (q) edge node[right]{}(r);

\end{tikzpicture}
          \caption{Graph $K_3 \cdot K_3$}
       \label{fig2}
\end{center}
\end{figure}
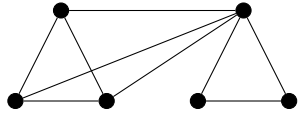

\noindent{\bf Remark:} This definition is a generalization of the graph $K K_n^2,$ introduced by Stevanovi\'c in \cite{steva}, where it is studied other spectral properties of this graph, as Laplacian energy.
The Figure \ref{fig2} shows the graph $K_3 \cdot K_3.$

\begin{Lem}
\label{lema2}
Let $G =  K_n \cdot K_n$ be a graph of order $2n.$
Then  the Laplacian spectrum of $ G$ is given by
$$ 0 ;  \hspace{0,2cm} 1 ;\hspace{0,2cm} n^{n-2}; \hspace{0,2cm} (n+1)^{n-1};\hspace{0,2cm} 2n.$$
\end{Lem}
{\bf Proof:}
Let  $G = K_n \cdot K_n$ be a graph of order $2n.$ Since that  $G$ can be viewed as the join $(K_{n} \cup K_{n-1}) \nabla K_1,$ the proof is similar to Lemma \ref{lema1}.

\begin{Thr}
\label{}
For each $n\geq 3,$  $G= K_n \cdot K_n$  is  $L$-borderenergetic and  $L$-noncospectral graph with $K_{2n}.$
\end{Thr}
{\bf Proof:}
Clearly $G= K_n \cdot K_n$ and $K_{2n}$ are $L$-noncospectral.
Let $\overline{d}$ be the average degree of $G.$ Since that  $\overline{d}= n$ then
 $E_{L}(G) = |2n -n| + (n-1) | n+1 -n| + (n-2) |n-n| +|1 -n| +|0-n| = 4n-2 = E_L(K_{2n}). $

\section{ Graphs with same normalized Laplacian energy}

In this section we present graphs which have the same normalized Laplacian energy $E_{\mathcal{L}}(G).$
Let be $K_n$ the complete graph on $n$ vertices.
For positive integer $b\geq 2$  we define the following classes of graphs:
\begin{itemize}
\item  The class of graphs  $    2K_2 \nabla  K_b.$
\item  The class of graphs  $    2K_2 \nabla  b K_1.$
\end{itemize}

We let $ \lambda_i^m$ denote the $m$- multiplicity of normalized laplacian eigenvalue $\lambda_i.$

\begin{Lem}
\label{lemma i}
Let be an integer positive $b\geq 2$ 
and $G = 2K_2 \nabla K_b$  a  graph of order $n=b+4.$ Then
$$ 0,    \frac{b}{b+1},  (\frac{b+2}{b+1})^2, (\frac{b+4}{b+3})^{b-1}, \frac{b}{b+1} + \frac{4}{b+3}.$$
are the $\mathcal{L}$- eigenvalues of $G.$ 
\end{Lem}
{\bf Proof:}
Let be $G_1 = 2K_2$ and $G_2= K_b$ the graphs of order $n=4$ and $m=b,$ respectively. Since  that $G_1$ is an $1$-regular graph and $G_2$ is an  $(b-1)$-regular graph, and the $\mathcal{L}$- eigenvalues of $G_1$ and $G_2$ are given by 
$\{ 0^2, 2^2\}$ and $\{0, (\frac{b}{b-1})^{b-1} \},$ respectively.
Taking $r=1$ and $s=b-1,$ then the result follows by Theorem (\ref{main4}).

\begin{Lem}
\label{lemma ii}
Let be an integer positive $b\geq 2$ 
and $G' = 2K_2 \nabla b K_1$  a  graph of order $n=b+4.$ Then
$$ 0,    \frac{b}{b+1},  (\frac{b+2}{b+1})^2, (1)^{b-1}, \frac{b}{b+1} + 1.$$
are the $\mathcal{L}$- eigenvalues of $G'.$ 
\end{Lem}
{\bf Proof:}
Similar to Lemma \ref{lemma i}.

\begin{Thr}
\label{main5}
Let be an integer positive $b\geq 2,$ 
 $G = 2K_2 \nabla K_b$  and $G' = 2K_2 \nabla bK_1$  graphs of order $n=b+4.$ Then $G$ and $G'$ are  $\mathcal{L}$-noncospectral and have the same normalized Laplacian energy. Furthermore
$$ E_{\mathcal{L}}(G) = E_{\mathcal{L}}(G')= \frac{2b+4}{b+1}.$$ 
\end{Thr}
{\bf Proof:}
Let be $G = 2K_2 \nabla K_b$ and $G'= 2K_2 \nabla bK_1$  graphs of order $n=b+4.$ Clearly $G$ and $G'$ are  $\mathcal{L}$-noncospectral.  Using that $E_{\mathcal{L}}(G)= \sum_{i=1}^n | \lambda_i - 1 |$ and Lemma \ref{lemma i}, follows 
$$E_{\mathcal{ L}}(G) = |0 -1| + | \frac{b}{b+1} -1| +2|\frac{b+2}{b+1} -1| +(b-1)| \frac{b+4}{b+3}-1| + | \frac{b}{b+1} + \frac{4}{b+3} -1|$$ 
$$ E_{\mathcal{L}}(G) = \frac{2b^2 +10b+12}{(b+1)(b+3)} = \frac{2(b+2)(b+3)}{(b+1)(b+3)} = \frac{2b +4}{b+1}.$$
If $G'= 2K_2 \nabla bK_1,$ by Lemma (\ref{lemma ii}), we have that
$$E_{\mathcal{ L}}(G') = |0 -1| + | \frac{b}{b+1} -1| +2|\frac{b+2}{b+1} -1| +(b-1) | 1-1| + | \frac{b}{b+1} + 1-1|$$ 
$$ E_{\mathcal{L}}(G') =    1+  \frac{3}{b+1} + \frac{b}{b+1}= \frac{2b +4}{b+1},$$
and then the result follows.

Now we present the general case.
Let be $K_n$ the complete graph on $n$ vertices.
For positive integers  $a,b\geq 2$  we define the following classes of graphs:
\begin{itemize}
\item  The class of graphs  $    aK_2 \nabla  K_b.$
\item  The class of graphs  $    aK_2 \nabla  b K_1.$
\end{itemize}

\begin{Lem}
\label{lemmazero}
Let be the integers positive $a, b\geq 2$ 
and $G = aK_2 \nabla K_b$  a  graph of order $n=2a +b.$ Then
$$ 0,    (\frac{b}{b+1})^{a-1},  (\frac{b+2}{b+1})^a, (\frac{2a+b}{2a +b-1})^{b-1}, \frac{b}{b+1} + \frac{2a}{2a +b -1}.$$
are the $\mathcal{L}$- eigenvalues of $G.$ 
\end{Lem}

\begin{Lem}
\label{lemmaum}
Let be the integers positive $a, b\geq 2$ 
and $G' = aK_2 \nabla bK_1$  a  graph of order $n=2a +b.$ Then
$$ 0,    (\frac{b}{b+1})^{a-1},  (\frac{b+2}{b+1})^a, (1)^{b-1}, \frac{b}{b+1} + 1.$$
are the $\mathcal{L}$- eigenvalues of $G'.$ 
\end{Lem}

\begin{Thr}
\label{main6}
Let be the integers positive $a,b\geq 2,$ 
 $G = aK_2 \nabla K_b$  and $G' = aK_2 \nabla bK_1$  graphs of order $n=2a +b.$ Then $G$ and $G'$ are  $\mathcal{L}$-noncospectral and have the same normalized laplacian energy. Furthermore
$$ E_{\mathcal{L}}(G) = E_{\mathcal{L}}(G')= \frac{2a+2b}{b+1}.$$ 
\end{Thr}
{\bf Proof:}
Let be $G = aK_2 \nabla K_b$ and $G'= aK_2 \nabla bK_1$  graphs of order $n=2a+b.$ Clearly $G$ and $G'$ are  $\mathcal{L}$-noncospectral.  Using that $E_{\mathcal{L}}(G)= \sum_{i=1}^n | \lambda_i - 1 |$ and Lemma \ref{lemmazero}, follows 
$$E_{\mathcal{ L}}(G) = |0 -1| +(a-1) |\frac{b}{b+1} -1| +a|\frac{b+2}{b+1} -1| +(b-1)| \frac{2a+b}{2a +b-1}-1| + | \frac{b}{b+1} + \frac{2a}{2a+b-1} -1|$$ 
$$ E_{\mathcal{L}}(G) = \frac{ 4a^2 +6ab + 2b^2 -2a  -2b}{(b+1)(2a +b-1)} = \frac{(2a+2b)( 2a + b-1)}{(b+1)(2a +b-1)} = \frac{2a +2b}{b+1}.$$
If $G'= aK_2 \nabla bK_1,$ by Lemma (\ref{lemmaum}), we have that
$$E_{\mathcal{ L}}(G') = |0 -1| + (a-1)| \frac{b}{b+1} -1| +a|\frac{b+2}{b+1} -1| +(b-1) | 1-1| + | \frac{b}{b+1} + 1-1|$$ 
$$ E_{\mathcal{L}}(G') =    1+  \frac{2a-1}{b+1} + \frac{b}{b+1}= \frac{2a +2b}{b+1},$$
and then the result follows.


\begin{thebibliography}{99}


\bibitem{Butler} 
S. Butler, Eigenvalues and Structures of Graphs, Ph.D. dissertation, University of California, San Diego, 2008.


\bibitem{Li2}
B. Deng, X. Li and I. Gutman,  More on  borderenergetic graphs, {\it Lin. Algebra Appl. \/}
{\bf 497} (2016) 199-208.




\bibitem{Li7}
 B. Deng, X. Li,    More  on $L$-borderenergetic  graphs, {\it MATCH Commun. Math. Comput. Chem. \/} {\bf 77} (2017) 115-127.



\bibitem{Godsil}
C. Godsil, G. Royle, {\it Algebraic Graph Theory\/}, Springer-Verlag, New York, 2001.




\bibitem{Gutman}
I. Gutman, B. Zhou,
Laplacian energy graph of a graph,
{\it Lin. Algebra  Appl. \/}, {\bf 414} (2006) 29--37.



\bibitem{JTT2015}
D. P. Jacobs, V. Trevisan, F. Tura, Eigenvalues and energy in
    threshold graphs, {\it Lin. Algebra Appl.\/} {\bf 465} (2015)
    412--425.





\bibitem{Merris2}
R. Merris,
 Laplacian  graph eigenvectors,
{\it Lin. Algebra  Appl \/}, {\bf 278} (1998) 221--236.



\bibitem{Gutman2015}
S. C. Gong,X. Li, G. H. Xu, I. Gutman, B. Furtula,
Borderenergetic graphs ,
{\it MATCH Commun. Math. Comput. Chem. \/}{\bf 74} (2015) 321-332.




\bibitem{Hou}
 Y. Hou, Q Tao, Borderenergetic threshold graphs, {\it MATCH Commun.
    Math. Comput. Chem.\/} {\bf 75} (2016) 253--262.



\bibitem{steva}
D. Stevanovi\'c, I. Stankovi\'c, M. Milosevi\'c, More on the relation between energy and   laplacian energy of graphs, {\it MATCH Commun. Math. Comput. Chem.\/}{\bf 61} n.2 (2009) 395--401.







\bibitem{Gutman2012}
X. Li, Y, Shi, I. Gutman, {\it Graph Energy\/}, Springer,
    New York, 2012.







\bibitem{Li}
X. Li, M. Wei, S. Gong, A computer search for the borderenergetic graphs of order 10,
{\it MATCH Commun. Math. Comput. Chem. \/}{\bf 74} (2015) 333-342.




\bibitem{Ramane}
H. S. Ramane, H. B. Walikar, Construction of equienergetic graphs,
{\it MATCH Commun. Math.  Comput. Chem. \/} {\bf 57}
(2007) 203--210.



\bibitem{Li3}
Z. Shao, F. Deng,  Correcting the number of borderenergetic graphs of order 10,
{\it MATCH Commun. Math. Comput. Chem. \/} {\bf 75} (2016) 263-266.



\bibitem{Li4}
 X. Li, H. Ma,   All hypoenergetic graphs with maximum degree at most 3, {\it Lin. Algebra Appl.\/} {\bf 431} (2009) 2127-2133.


\bibitem{Li5}
 X. Li, H. Ma,   All  connected graphs with maximum degree at most 3 whose energies are equal to
 the number of vertices, {\it MATCH Commun. Math. Comput. Chem. \/} {\bf 64(1)}  (2010) 7-24.


\bibitem{Li6}
 X. Li, H. Ma,   Hypoenergetic and strongly hypoenergetic k-cyclic graphs, {\it MATCH Commun. Math. Comput. Chem. \/} {\bf 64(1)} (2010) 41-60.




\bibitem{tura}
F. Tura, $L$-Borderenergetic graphs,
{\it MATCH Commun. Math.  Comput. Chem. \/} {\bf 77}
(2017) 37--44.


\end{thebibliography}
\end{document}